\documentclass[12pt]{amsart}
\usepackage{amsmath,amsthm,amsfonts,amssymb,times}
\DeclareMathAlphabet{\curly}{U}{rsfs}{m}{n}

\theoremstyle{remark}

\theoremstyle{plain}

\newtheorem{lem}{Lemma}[section]
\newtheorem{thm}{Theorem}
\newtheorem{cor}{Corollary}
\numberwithin{equation}{section}

\newcommand{\R}{{\mathbb R}}   
\newcommand{\ZZ}{{\mathbb Z}}  

\newcommand{\EE}{\ensuremath{\mathbf{E}}}  
\newcommand{\PP}{\ensuremath{\mathbf{P}}}  
\newcommand{\LL}{\curly L}

\newcommand{\be}{\begin{equation}}
\newcommand{\ee}{\end{equation}}
\newcommand{\benn}{\begin{equation*}}
\newcommand{\eenn}{\end{equation*}}

\newcommand{\lam}{\ensuremath{\lambda}}

\renewcommand{\a}{\ensuremath{\alpha}}
\renewcommand{\b}{\ensuremath{\beta}}
\newcommand{\g}{\ensuremath{\gamma}}

\renewcommand{\(}{\left(}
\renewcommand{\)}{\right)}

\newcommand{\pfrac}[2]{\left(\frac{#1}{#2}\right)}

\newcommand{\RR}{\ensuremath{\widetilde{R}}}

\newcommand{\tS}{\ensuremath{\widetilde{S}}}
\newcommand{\tT}{\ensuremath{\widetilde{T}}}

\begin{document}

\title{Sharp probability estimates for random
walks with barriers}
\author{Kevin Ford}
\address{
Department of Mathematics \\
University of Illinois  \\
Urbana, IL 61801}
\email{ford@math.uiuc.edu}
\date{11 July 2008}

\thanks{The author was supported by NSF grants DMS-0301083 and DMS-0555367.}



\begin{abstract}
We give sharp, uniform estimates for the probability that a random
walk of $n$ steps on the reals avoids a half-line $[y,\infty)$ given 
that it ends at the point $x$. The estimates hold for general
continuous or lattice distributions provided the 4th moment is finite.
\end{abstract}

\maketitle

%
\section{Introduction}
%

Let $X_1, X_2, \ldots$ be independent, identically distributed
random variables with mean $\EE X_1 = 0$ and variance $\EE
X_1^2 = 1$.  Let $S_0=T_0=0$ and for $n\ge 1$ define
$$
S_n = X_1 + \cdots + X_n
$$
and
$$
T_n = \max (0,S_1,\ldots,S_n).
$$
The estimation of the distribution of $S_n$ for general random
variables has a long and rich history (see e.g. \cite{GK}).


The distribution of $T_n$ was found more recently.
In 1946, Erd\H os and Kac \cite{EK} showed that
$$
\lim_{n\to \infty} \PP[T_n \le x\sqrt{n}] = 
2\Phi(x)-1
$$
uniformly in $x\ge 0$, where
$$
\Phi(x) = \frac{1}{\sqrt{2\pi}} \int_{-\infty}^x e^{-t^2/2}\, dt
$$
is the distribution function for the normal distribution.
Several estimates for the error term 
have been proved on the assumption that $\EE |X_1|^3 < \infty$, 
the best uniform bound (and best possible uniform
bound) being the result of Nagaev \cite{N}
$$
\PP[T_n \le x\sqrt{n}] = 2\Phi(x) - 1 + O(1/\sqrt{n}),
$$
uniformly in $x\ge 0$ (the constant implied by the $O-$symbol depends only
on $\EE |X_1|^3$).  Sharper error terms are possible when $|x|\ge
1$, see e.g. Arak \cite{A} and Chapter 4 of \cite{AN}.

We are interested here in approximations of the conditional probability
$$
R_n(x,y) = \PP[T_{n-1} < y | S_n = x]
$$
which are sharp for a wide range of $x,y$.
By the invariance principle, we expect
$$
R_n(u\sqrt{n},v\sqrt{n}) \to 1 - e^{-2v(v-u)} \qquad (n\to \infty)
$$
for $u,v$ fixed, $u\le v$ and $v\ge 0$, since this holds for the case of
Bernoulli random variables (see \eqref{Bern} below).

Before stating our results, we motivate the study of $R_n(x,y)$
with three examples, two of which are connected with empirical processes.

%
\section{Three examples}\label{sec:examples}
%

The example which is easiest to analyze is the case of a 
simple random walk with Bernoulli steps.
Let $X_1, X_2, \ldots$ satisfy $\PP[X_i=1]=\PP[X_i=-1]=\frac12$.  
By the reflection principle, the number of walks of $n$ steps
for which $T_n \ge y$ and $S_n=x$ is equal to the number of walks of $n$ steps
with $S_n=2y-x$ (by inverting $X_{k+1},\ldots,X_n$, where $k$ is
the smallest index with $S_k=y$).  Thus, if $n$ and $x$ have the same
parity, then
$$
R_n(x,y) = 1 -
\frac{\binom{n}{\frac{n+x-2y}{2}}}{\binom{n}{\frac{n-x}{2}}}. 
$$
This includes as a special case a version of the classical Bertrand
ballot theorem from 1887.  Two candidates in an election
receive $p$ and $q$ votes, respectively, with $p\ge q$.  If the votes
are counted in random order, the probability that the first candidate
never trails in the counting is
$$
R_{p+q}(q-p,1) = \frac{p-q+1}{p+1}.
$$

More generally, suppose $1\le y \le n/2$, $-n/2 \le x < y$ and $2y-x\le n/2$.
Writing $\b=(2y-x)/n$ and $\a=x/n$, so that $\b>\a>0$, we obtain by
Stirling's formula,
\begin{align*}
R_n&(x,y) = 1 -
  \frac{\binom{n}{\frac{n}{2}(1+\b)}}{\binom{n}{\frac{n}{2}
  (1+\a)}} \\
&= 1 - (1+O(1/n)) \sqrt{\frac{1-\a^2}{1-\b^2}} \( \frac{(1+\a)^{1+\a}
    (1-\a)^{1-\a}}{(1+\b)^{1+\b}(1-\b)^{1-\b}} \)^{n/2} \\
&=1 - (1+O(1/n)) \sqrt{\frac{1-\a^2}{1-\b^2}} \exp \left\{ \frac{n}{2}
\( \a^2 - \b^2 +O\(\a^4+\b^4\) \)  \right\}.
\end{align*}
If $x=O(\sqrt{n})$ and $y-x=O(\sqrt{n})$, then $\a = O(n^{-1/2})$ and
$\b=O(n^{-1/2})$ and we have
\be\label{Bern}
R_n(x,y) = 1 - (1+O(1/n)) \exp\{ \tfrac{n}{2}(\a^2-\b^2) \}
= 1 - e^{-2y(y-x)/n} + O(1/n).
\ee

Two special cases are connected with empirical processes.
Let $U_1,\ldots,U_n$ be independent random
variables with uniform distribution in $[0,1]$, suppose
$F_n(t) = \frac{1}{n} \sum_{U_i\le t} 1$ is their
empirical distribution function and
$0 \le \xi_1 \le \cdots \le \xi_n \le 1$ are their order statistics.

In his seminal 1933 paper
\cite{Kol33} on the distribution of the statistic
$$
D_n = \sqrt{n} \sup_{0\le t\le 1} |F_n(t)-t|,
$$
Kolmogorov related the problem to a similar conditional probability
for a random walk.
Specifically, let $X_1,X_2,\ldots,X_n$ be independent random
variables with discrete distribution 
\be\label{KolXj}
\PP [ X_j=r-1 ] = \frac{e^{-1}}{r!} \qquad (r=0,1,2,\ldots)
\ee
Kolmogorov proved that for integers $u\ge 1$,
\begin{align*}
\PP ( \sup_{0\le t\le 1} |F_n(t)-t| \le u/n ) 
&= \frac{n! e^n}{n^n} \PP \( \max_{0\le
  j\le n-1} |S_j| < u, S_n=0 \) \\
&= \PP \(  \max_{0\le j\le n-1} |S_j| < u \, \Big| S_n=0 \).
\end{align*}

Consider next
$$
Q_n(u,v) = \PP[ \xi_i \ge \tfrac{i-u}{v} \; (1\le i\le n)] = 
\PP\( F_n(t) \le \frac{vt+u}{n} \;\; (0\le t\le 1) \) 
$$
for $u\ge 0, v>0$.
Smirnov in 1939 proved the asymptotic $Q_n(\lam\sqrt{n},n)\to 1-e^{-2\lam^2}$
as $n\to\infty$ for fixed $\lambda$.
Small modifications to Kolmogorov's proof yield, 
for {\it integers} $u\ge 1$ and for $n\ge 2$, that
$$
Q_n(u,n) =  R_n(0,u)
$$
for the variables $X_j$ given by \eqref{KolXj}.
When $v\ne n$, however, it does not seem possible to express $Q_n(u,v)$
in terms of these variables $X_j$.

In \cite{F}, new bounds on $Q_n(u,v)$ were proved
and applied to a problem of the
distribution of divisors of integers (see also articles \cite{kol1},
\cite{kol2} for more about this application).
A more precise uniform estimate was proved in \cite{smir}, namely 
\be\label{Qnuv}
Q_n(u,v) = 1-e^{-2uw/n} + O\pfrac{u+w}{n} \qquad (n\ge 1, u\ge 0, w\ge 0),
\ee
where $w=u+v-n$ and the constant implied by the $O-$symbol is
independent of $u,v$ and $n$.  This was accomplished using
$X_j=1-Y_j$, where
$Y_1, Y_2,\ldots$ are independent random variables with
exponential distribution, i.e. with density function 
$f(x) = e^{-x}$ for $x\ge 0$, $f(x)=0$ for $x<0$.
Letting $W_k=Y_1+\cdots+Y_{k}$, R\'enyi \cite{Re} whowed that
$$
(\xi_1,\xi_2,\cdots,\xi_n) \text{ and }
\( \frac{W_1}{W_{n+1}}, \frac{W_2}{W_{n+1}},
 \cdots, \frac{W_n}{W_{n+1}} \)
$$
have the same distribution.  An easy consequence is
$$
Q_n(u,v) = \PP \bigl[W_j-j\ge -u\; (1\le j\le n)\; | \; W_{n+1} = v\bigr] 
= R_{n+1}(n+1-v,u).
$$

%
\section{Statement of the main results}
%

Our aim in this paper is to prove a result analogous to
\eqref{Bern} and \eqref{Qnuv}
for sums of very general random variables $X_1$.
We will restrict ourselves to random variables with either a continuous
or lattice distribution, to maintain control of the density function
of $S_n$.  Let $F$ be the distribution
function of $X_1$ and let $F_n$ the distribution function of $S_n$ for
$n\ge 1$.  Let $\phi(t)=\EE e^{it X_1}$ be the characteristic function
of $X_1$.

We henceforth assume that either\benn\label{hypC}
\begin{cases} X_1 \text{ has a continuous distribution and} &\\
\exists n_0 : \, \int |\phi(t)|^{n_0}\, dt < \infty &\end{cases} \tag{C}
\eenn
or that
\benn\label{hypL}
X_1 \text{ has a lattice distribution.} \tag{L}
\eenn
If (L), let $f(x)=\PP(X_1=x)$, $f_n(x)=\PP(S_n=x)$ and $n_0=1$.
We also suppose the support of $f$ is contained in
the lattice $\LL = \{ \g+m\lam : m\in \ZZ \}$, where $\lam$ is the maximal
span of the distribution (the support of $f$ is not contained in any
lattice $\{\g'+m\lam':m\in \ZZ\}$ with $\lam'>\lam$).  The support of $f_n$ is
then contained in the lattice $\LL_n = \{ n\g+m\lam:m\in\ZZ\}$.
If (C), let $f$ be the
density function of $X_1$, let $f_n$ the density function of $S_n$,
define $\LL=\R$ and $\LL_n=\R$.

Define the moments
$$
\a_u = \EE X_1^u, \qquad \b_u = \EE |X_1|^u.
$$

%

In what follows, the notation $f=O(g)$ for functions $f,g$ means that
for some constant $c>0$, $|f| \le c g$ for all values of the domain of
$f$, which will usually be given explicitly.
Unless otherwise specified, $c$ may depend only 
on the distribution of $X_1$, but not on any other parameter.
Sometimes we use the Vinogradov notation $f \ll g$ which means
$f=O(g)$.  As $R_n(x,y)$ is only defined when $f_n(x)>0$,
when $f_n(x)=0$ we define $R_n(x,y)=1$.

%
%

\begin{thm}\label{posz}
Assume (C) or (L), $\b_u<\infty$ for some $u>3$,
and let $M>0$.  
Uniformly in $n\ge 1$, $0\le y \le M \sqrt{n}$, $0\le z\le M \sqrt{n}$
with $y\in\LL_n$, $y-z\in\LL_n$ and $f_n(y-z)>0$, 
$$
R_n(y-z,y) = 1 - e^{-2yz/n} + O\( \frac{y+z+1}{n} +
\frac{1}{n^{\frac{u-2}{2}}}\).
$$
Here the constant implied by the $O-$symbol depends on the
distribution of $X_1$, $u$ and also on $M$, but not on $n,y$ or $z$.
\end{thm}

\begin{cor}
Assume (C) or (L) and $\beta_u<\infty$ for some $u>3$.  
For $w\le v$ and $v\ge 0$, 
$$
R_n(w\sqrt{n}, v\sqrt{n}) = 1 - e^{-2v(v-w)} + O(n^{-1/2}),
$$
the constant implied by the $O-$symbol depending on $\max(v,v-w)$
and on the distribution of $X_1$.
\end{cor}

\begin{cor}
Assume (C) or (L), and $\beta_4<\infty$.  If $y$ and $z$ satisfy
$y\to \infty$, $y=o(\sqrt{n})$, $z\to \infty$,  and
$z=o(\sqrt{n})$ as $n\to \infty$, then
$$
\lim_{n\to \infty} \frac{R_n(y-z,y)}{2yz/n} = 1.
$$
\end{cor}

All three examples given in section \ref{sec:examples} staisfy the hypotheses of
Theorem 1 and the two corollaries.  Indeed, for these examples all
moments of $X_1$ exist. 

Using ``almost sure invariance'' principles or ``strong approximation''
theorems (see e.g. \cite{CR}, \cite{Ph}), one can approximate the
walk $(S_n)_{n\ge 0}$ with a Wiener process $W(n)$.   Assuming
that $\b_4<\infty$ and no higher moments exist, one has
$S_n-W(n)=o(n^{1/4})$ almost surely, the exponent $1/4$ being best
possible (cf. \cite{CR}, Theorems 2.6.3, 2.6.4).  This rate of approximation
is, however, far too weak to prove results as strong as Theorem
\ref{posz}.

In section \ref{sec:basic_estimates}, we list some required 
estimates for $f_n(x)$.   Section
\ref{sec:recurrence} contains two recursion formulas for $R_n(x,y)$.
Although our main interest is in the case when $y\ge x$, we shall need
estimates when $y<x$ in order to prove Theorem \ref{posz}.  
This is accomplished in \S \ref{sec:negz}.  Finally, in \S \ref{pfthm}, we
prove Theorem \ref{posz}.
It is critical to our analysis that the
densities $f_n(x)$ have regular behavior, and the hypotheses (C), (L)
and $\b_u<\infty$
ensures that this is the case for $|x| = O(\sqrt{n})$.  Extending the
range of validity of the asymptotic for $R_n(x,y)$ beyond the range of $(x,y)$
covered by Theorem \ref{posz} would require that we have more precise estimates
for $f_n(x)$ for $|x|$ of larger order than $\sqrt{n}$.  In specific
cases, such as the exponential distribution, normal distribution or
binomial distribution, exact expressions for $f_n(x)$
make it possible to achieve this goal (see e.g. \eqref{Qnuv}).

It is of some interest to know if $\beta_4<\infty$ really is
a necessary condition for Theorem 1 to hold.  Recently, Addario-Berry
and Reed \cite{ABR} showed (as a special case of their Theorem 1), 
for an arbitrary lattice random variable $X_1$
with zero mean and finite variance, that
\be\label{ABR}
\frac{yz}{n} \ll R_n(y-z,y) \ll \frac{yz}{n} \qquad (1\le y,z\le \sqrt{n},
n\ge n_0),
\ee
the constants implied by the $\ll$-symbols and $n_0$ depending on the
distribution of $X_1$.  
The same
proof gives \eqref{ABR} under hypotheses (C) and $\beta_2 < \infty$;
see \eqref{roughsmally} below
(for non-lattice variables, the authors prove analogous bounds for
the probability that $T_n<y$ given that $y-z-c \le S_n\le y-z$, for a
fixed $c>0$).   When $y=1$, the upper bound in \eqref{ABR} is the same as
the conclusion as Theorem \ref{posz},
but is proved under a weaker hypothesis.  When $y$ is larger,
however, the error term in the conclusion of Theorem \ref{posz} can be
of much lower order than the main term, and a hypothesis stronger
than $\beta_2<\infty$ should be required.
Addario-Berry and Reed also construct examples of variables $X_1$ where
 $\EE X_1^2=\infty$ or $k/\sqrt{n} \to \infty$, while $R_n(-k,1)$
is not of order $k/n$.

%
\section{Estimates for density functions}\label{sec:basic_estimates}
%

At the core of our arguments are approximations of the density
function $f_n(x)$.  This is the only part of the proof which uses the
hypothesis on $\phi(t)$ from (C).

\begin{lem}\label{fn1} Assumer (C) or (L), and $\beta_2=1$.  Then,
  uniformly for $n\ge n_0$ and all $x$,
\be\label{fnbeta2}
f_n(x) \ll \frac{1}{\sqrt{n}}.
\ee
Assume $3\le u\le 4$, $\b_u<\infty$, and (C) or (L).
Then, uniformly for $n\ge n_0$ and $x\in\LL_n$,
\begin{align*}
f_n(x) &= \frac{e^{-x^2/2n} }{\sqrt{2\pi n}}
 \left[ 1 + O\( \frac{|x|}{n} + \frac{|x|^3}{n^2} \) \right] +
 O(n^{(1-u)/2}) \\
&= \frac{ e^{-x^2/2n}}{\sqrt{2\pi n}} +
O \( \frac{|x|}{n^{3/2}(1+x^2/n)} + n^{(1-u)/2} \).
\end{align*}
\end{lem}

\begin{proof}
We apply results from \cite{GK}, \S 46, \S 47 and \S 51.
Assume (C).  By the proof of Theorem 1 in \S 46, 
we may replace conditions 1), 2) of \S 46, Theorem 1 and
the theorem in \S 47
with the hypothesis that $n_0$ exists.  Note that these theorems are
only stated with the hypothesis that $\beta_u$ exists for intergal
$u$, but straightforward modification of the proofs yields the above
inequalities for real $u\in [3,4]$: Start with the inequality
$e^{it} = 1 + it - \frac12 t^2 - \frac{i}{6} t^3 + O(|t|^u),$
which follows from Taylor's formula for $|t|\le 1$ and the triangle
inequality for $|t|>1$.  Consequently,
$$
\phi(t) = 1 - \frac12 t^2 - \frac{i \a_3}{6} t^3 + O(|t|^u)
$$
and hence, for $|t|$ small enough,
\begin{align*}
\phi^n(t) &= \exp \left[ - \frac{nt^2}{2} - \frac{i\a_3n}{6} t^3 + O(n
  |t|^u) \right] \\
&= e^{-nt^2/2} \left[ 1 - \frac{i \a_3 n}{6} t^3 + O\( t^6 n^2 e^{O(|t|^3n)}+
  |t|^u n e^{O(|t|^un)} \) \right].
\end{align*}
Here we used the inequalities $|e^v-1| \le |v| e^{|v|}$ and
$|e^v-1-v| \le |v|^2 e^{|v|}$.  Therefore,
\be\label{phint}
\left|\phi^n(t) - e^{-nt^2/2} \( 1 - \frac{i\a_3 n}{6} t^3 \) \right|
\ll (t^6 n^2 + |t|^u n) e^{-nt^2/4} \quad (|t| \le c)
\ee
for some $c>0$.
In the proofs in \S 46, \S 47 and \S 51, use \eqref{phint} in place of
Theorem 1 of \S 41.
\end{proof}


%
\section{Recursion formulas}\label{sec:recurrence}
%

It is convenient to work with the density function
$$
\RR_n(x,y) = f_n(x) R_n(x,y) = \PP[T_{n-1} < y, S_n=x].
$$
The last expression stands for $\frac{d}{dx} \PP[T_{n-1} < y, S_n\le
  x]$ when (C) holds.
Notice that if $f_n(x)=0$, then $\RR_n(x,y)=0$ by our convention.

\begin{lem}\label{recur1}
Assume (C).  Then, for $n\ge 2$, $y\ge 0$ and $s\ge 0$,
$$
\RR_n(y+s,y) = \int_0^\infty f(s+t) \RR_{n-1}(y-t,y)\, dt.
$$ 
If (L), then for $n\ge 2$, $y>0$, $s\ge 0$ and $y+s\in \LL_n$,
$$
\RR_n(y+s,y) = \sum_{\substack{s+t\in\LL_n \\ t>0}} f(s+t) \RR_{n-1}(y-t,y).
$$
\end{lem}

\begin{proof}
If $S_n=y+s$ and $T_{n-1} < y$, then $X_n=s+t$ where $t > 0$.
\end{proof}

Lemma \ref{recur1} expresses $\RR_n(x,y)$ with $x\ge y$ in terms of
$\RR_{n-1}(x,y)$ with
$x \le y$.  The next lemma works the other direction,
and is motivated by the reflection principle:
a walk that crosses the point $y$ and ends up at
$S_n=x$ should be about as likely as a walk that ends up at $S_n=2y-x$
(by inverting the part of the walk past the first crossing of $y$).
We thus expect that for $x<y$, 
$$
\RR_n(x,y) \approx f_n(x) - f_n(2y-x).
$$

\begin{lem}\label{recur2}  Assume $n\ge 2$, 
$y > 0$ and $a\ge 0$.  If (C), then for any $x$
\begin{align*}
\RR_n(x,y) = f_n(x) & - f_n(y+a) + \RR_n(y+a,y) \\
& + \int_0^\infty \sum_{k=1}^{n-1} \RR_k(y+\xi,y) \(
f_{n-k}(a-\xi)-f_{n-k}(x-y-\xi) \)\, d\xi.
\end{align*}
If (L), then for $x,y+a\in \LL_n$,
\begin{align*}
\RR_n(x,y) = f_n(x) & - f_n(y+a) + \RR_n(y+a,y) \\
& + \sum_{k=1}^{n-1} \sum_{\substack{y+\xi\in\LL_k \\ \xi\ge 0}}
 \RR_k(y+\xi,y) \( f_{n-k}(a-\xi)-f_{n-k}(x-y-\xi) \).
\end{align*}
\end{lem}

\begin{proof}  First, we have
\begin{align*}
\RR_n(x,y) &= f_n(x) - \PP[T_{n-1} \ge y, S_n=x] \\
&= f_n(x) - f_n(y+a) + f_n(y+a) -  \PP[T_{n-1} \ge y, S_n=x].
\end{align*}
If $T_j \ge y$, then there is a unique $k$, $1\le k\le j$, for which
$T_{k-1} < y$ and $S_{k}\ge y$.  Thus,
$$
f_n(y+a) = \sum_{k=1}^n \PP[T_{k-1} < y, S_k \ge y, S_n=y+a].
$$
If (C) then
\begin{align*}
f_n(y+a)
&= \sum_{k=1}^{n-1} \int_{0}^\infty \PP[T_{k-1} < y, S_k = y+\xi,
    S_n=y+a] \, d\xi + \PP[T_{n-1} < y, S_n=y+a] \\
&= \sum_{k=1}^{n-1}  \int_{0}^\infty \RR_{k}(y+\xi,y) f_{n-k}(a-\xi)\,
  d\xi + \RR_n(y+a,y).
\end{align*}
Likewise, if (L) then
$$
f_n(y+a) = \RR_n(y+a,y) + 
\sum_{k=1}^{n-1}  \sum_{\substack{y+\xi\in\LL_k \\ \xi\ge 0}}
\RR_{k}(y+\xi,y) f_{n-k}(a-\xi).
$$
In the same way
\begin{align*}
\PP[T_{n-1}\ge y, S_n=x] &= \sum_{k=1}^{n-1} \PP[T_{k-1} < y, S_k \ge y, S_n=x) \\
&= \sum_{k=1}^{n-1} \begin{cases} 
\int_{0}^\infty \RR_{k}(y+\xi,y) f_{n-k}(x-y-\xi)\,  d\xi & \text{ if
  (C)}\\
\sum_{\substack{y+\xi\in\LL_k \\ \xi\ge 0}} \RR_{k}(y+\xi,y) 
f_{n-k}(x-y-\xi) & \text{ if (L).}
\end{cases}
\end{align*}
\end{proof}

Motivated by the reflection principle, we will apply Lemma
\ref{recur2} with $a$ close to $y-x$.   The integral/sum over
$\xi$ is then expected to be small, since
$f_{n-k}(y-x-\xi)-f_{n-k}(x-y-\xi)$ should be small when $\xi$
is small (by Lemma \ref{fn1}) and
$\RR_k(y+\xi,y)$ should be small when $\xi$ is large.
This last fact is crucial to our argument, and we develop 
the necessary bounds in the next section.

%
\section{Rough Estimates}\label{sec:negz}
%

Roughly speaking, if $T_{n-1}<y$ and $S_n=y+s$ with $s\ge 0$, then $S_{n-1}$
will be close to $y$ with high probability.  The probability that $T_{n-1}<y$
is about $\min(1,y/\sqrt{n})$ (see Lemma \ref{PPlem} below)
On the other hand, if
$y/\sqrt{n}$ is large, then $S_{n-1}\approx y$ is a rare event.
Therefore, as a function of $y$, $\RR_n(y+s,y)$ should
increase linearly in $y$ for $1\le y\le \sqrt{n}$, then decrease for
larger $y$. 

We begin with a lemma concerning the distribution of $T_n$.
Part (1) is due to Kozlov (\cite{Koz}, Theorem A, (13)) and part (2)
was proved by Pemantle and Peres (\cite{PP}, Lemma 3.3 (ii)).  In
\cite{PP}, (2) is stated only for $h=0$, but
the same proof gives the more general inequality.

\begin{lem}\label{PPlem}
Assume $X_1$ is any random variable with $\beta_2=1$.  Then
\begin{enumerate}
\item $\PP \{ T_n\le h \} \ll (h+1)/\sqrt{n}$.
\item $\EE \{ S_n^2 | T_n \le h \} \ll n$, uniformly in $h\ge 0$.
\end{enumerate}
\end{lem}

\begin{thm}\label{thm2}
Assume (C) or (L), $\beta_2=1$ and $n\ge 1$.  For all $y\ge 0$, $z\ge 0$,
we have
\be\label{roughsmally}
\RR_n(y-z,y) \ll \frac{\min(y+1,\sqrt{n}) \min(z+1,\sqrt{n})}{n^{3/2}}.\tag{a}
\ee
If $n\ge 3n_0$, $y\ge \sqrt{n}$ and $0\le z \le y/2$, then
\be\label{roughlargey}
\RR_n(y-z,y) \ll \frac{\min(z+1,\sqrt{n})}{y^2}.\tag{b}
\ee
\end{thm}

\begin{proof}
The proof of \eqref{roughsmally} follows the upper bound proof of
Theorem 1 from \cite{ABR}.  The idea is to consider simultaneously the
random walk $0, S_1, S_2, \ldots$ and the ``reverse'' walk $0,\tS_1,
\tS_2, \ldots$, where $\tS_k = - (X_n + X_{n-1} + \cdots +
X_{n-k+1})$.  Let $\tT_n = \max(0,\tS_1,\ldots,\tS_n)$.  Note that
$T_n\le y$ and $S_n=y-z$ imply $\tT_n \le z$.  

Inequality \eqref{roughsmally} is trivial for $1\le n < 3n_0$.
Let $n\ge 3n_0$, put $a=\lfloor n/3 \rfloor$ and 
$b=n-a$.  Then $\RR_n(y-z,y) \le
\PP(E_1, E_2, E_3)$, where $E_1 = \{ T_a \le y \}$, $E_2 = \{ \tT_a
\le z\}$ and $E_3= \{ S_n = y-z \}$.  Think of the random walk $0, S_1,
\ldots, S_n$ as the union of three independent subwalks: one
consisting of the first $a$ steps, one consisting of steps numbered 
$a+1$ to $b$, 
and one consiting of the last $a$ steps reversed.
Note that $E_3 = \{S_b - S_a = y - z - S_a + \tS_a\}$.
Since $S_b-S_a$ is independent of $S_a$, $\tS_a$ and of events $E_1$
and $E_2$, we have by \eqref{fnbeta2}
$$
\PP (E_3 | E_1, E_2) \le \sup_w f_{b-a}(w) \ll n^{-1/2}.
$$
As $E_1$ and $E_2$ are independent, we have by Lemma \ref{PPlem} part (1)
$$
\RR_n(y-z,y) \le \PP E_1 \, \PP E_2 \, \PP \{ E_3 | E_1, E_2 \}
\ll \frac{\min(y+1,\sqrt{n}) \min(z+1,\sqrt{n})}{n^{3/2}}.
$$

To prove \eqref{roughlargey}, we observe that $S_n=y-z \ge y/2$.
Thus, $S_a \ge y/6$, $S_b-S_a \ge y/6$ or $\tS_a \le -y/6$.
Suppose first that $S_a \ge y/6$.  Replace $E_1$ by
$E_1' = \{ S_a \ge y/6 \}$ in the above argument and note that
$$
\PP \{ S_a \ge y/6 \} \le \frac{\EE S_a^2}{(y/6)^2} = \frac{36a}{y^2}
\ll \frac{n}{y^2}.
$$
Arguing as in the proof of \eqref{roughsmally}, we find that
$$
\PP \{ T_n \le y, S_n=y-z, S_a \ge y/6\} \ll \pfrac{n}{y^2}
\frac{\min(z+1,\sqrt{n})}{n} \ll \frac{\min(z+1,\sqrt{n})}{y^2}.
$$
Next, suppose that $S_b-S_a \ge y/6$.  In the above argument, replace
$E_1$ with $E_1''=\{ S_b - S_a \ge y/6\}$.  Then $E_3 = \{
S_a = y-z-(S_b-S_a)+\tS_a\}$.
Again, $\PP \{ E_3|E_1'',E_2\} \le \sup_w f_a(w) \ll n^{-1/2}$ and we
obtain
$$
\PP \{ T_n \le y, S_n=y-z, S_b-S_a \ge y/6\} \ll
\frac{\min(z+1,\sqrt{n})}{y^2}.
$$
Finally, suppose $\tS_a \le -y/6$.  Replace $E_2$ with $E_2'= \{ \tS_a
\le -y/6, \tT_a \le z\}$.  Here we use the trivial bound $\PP E_1 \le
1$ and deduce
$$
\PP \{ T_n \le y, S_n=y-z, \tS_a \le -y/6 \} \le \PP E_1 \; \PP E_2'
\; \PP \{ E_3 | E_1, E_2' \} \ll n^{-1/2} \PP E_2'.
$$
By Markov's inequality and Lemma \ref{PPlem} parts (1) and (2), 
$$
\PP E_2' \le \PP \{ \tT_a\le z\} \PP \{ \tS_a\ge y/6 | \tT_a\le z \} \le
\PP \{ \tT_a \le z \} \frac{\EE \{ \tS_a^2 | \tT_a \le z
  \}}{(y/6)^2} 
\ll\frac{\min(z+1,\sqrt{n})}{\sqrt{n}} \, \cdot \, \frac{n}{y^2}.
$$
This completes the proof of \eqref{roughlargey}.
\end{proof}

Combining Theorem \ref{thm2} with Lemma \ref{recur1}
gives us useful bounds on $\RR_n(y+\xi,y)$ when $\xi \ge 0$.

\begin{thm}\label{xlarge}
Assume (C) or (L), and $\beta_2=1$.
Suppose $y \ge 0$ and $\xi\ge 0$.  Then
$$
\RR_n(y+\xi,y) \ll \frac{y+1}{n^{3/2}}\int_{0}^\infty (t+1) f(\xi+t)\,
dt \ll \frac{y+1}{n^{3/2}}.
$$
If $n\ge 3n_0+1$ and $y>\sqrt{n}$, then
$$
\RR_n(y+\xi,y) \ll \frac{1}{y^2}  
\int_{0}^\infty (t+1) f(\xi+t)\, dt + \frac{1-F(\xi+y/2)}{n^{1/2}}.
$$
\end{thm}

\begin{proof}
Apply Lemma \ref{recur1} and Theorem \ref{thm2} \eqref{roughsmally} for the
first part, and observe that the integral is $\le \EE |X_1|$.
For the second part, use  Theorem \ref{thm2}
\eqref{roughlargey} for $t\le y/2$, and
 $\RR_{n-1}(y-t,y) \ll n^{-1/2}$ for $t>y/2$.  
\end{proof}

%
\section{Proof of Theorem \ref{posz}}\label{pfthm}
%

We begin by proving a lemma which is of independent interest.

\begin{lem}\label{jumpsize}
Assume $\b_u<\infty$ for some $u\ge 2$, and $y \ge 0$.  If (C) then
$$
\sum_{n=1}^\infty \int_0^\infty \xi^{u-2} \RR_n(y+\xi,y)\, d\xi = O(1).
$$
If (L) then
$$
\sum_{n=1}^\infty \sum_{\substack{y+\xi\in\LL_n \\ \xi\ge 0}} 
\xi^{u-2} \RR_n(y+\xi,y) = O(1).
$$
\end{lem}

\begin{proof}  Assume (C).  First,
$$
\sum_{n=1}^{3n_0} \int_0^\infty \xi^{u-2} \RR_n(y+\xi,y)\, d\xi \le
\sum_{n=1}^{3n_0} \int_0^\infty \xi^{u-2} f_n(y+\xi)\, d\xi \ll
\sum_{n=1}^{3n_0} \EE |S_n|^{u-1} \ll 1.
$$
By Theorem \ref{xlarge},
\begin{align*}
\sum_{n\ge 3n_0+1} \int_0^\infty \xi^{u-2}  \RR_n(y+\xi,y)\, d\xi &\ll 
\( \sum_{3n_0+1\le n\le y^2+1} \frac{1}{y^2} + \sum_{n>y^2+1}
\frac{y+1}{n^{3/2}}\) \\
&\qquad \qquad \times  \int_0^\infty (t+1) \int_0^\infty
\xi^{u-2} f(\xi+t)\, d\xi\, dt \\
&\qquad + \sum_{n\le y^2+1} \frac{1}{n^{1/2}} \int_0^\infty \xi^{u-2}
\int_{\xi}^\infty f(v+y/2)\, dv \, d\xi \\
&\ll \EE (|X_1|^u+|X_1|^{u-1}) + (y+1) \int_0^\infty v^{u-1}
f(v+y/2)\, dv \\
&\ll 1 + \EE |X_1|^u \ll 1.
\end{align*}
The proof when (L) holds is similar.
\end{proof}

\noindent
{\bf Remark.}
A random walk $S_0, S_1,
\ldots$ with $\beta_2=1$ crosses the point $y$ with probability 1.
There is a unique $n$ for which
$T_{n-1} < y$ and $S_n \ge y$, and   Lemma \ref{jumpsize} states that
$\EE (S_n-y)^{u-2}= O(1)$.

We now prove Theorem \ref{posz} (again showing the details only for
the case of (C) holding).
It suffices to assume that $n$ is sufficiently large.
Let $n \ge 10n_0$ and put
$x=y-z$.   By Lemma \ref{recur2} with $a=z$,
\be\label{Rmain}
\begin{split}
\RR_n(x,y)&=f_n(x) - f_n(y+z) + \RR_n(y+z,y) \\
&\qquad\qquad + \int_0^\infty \sum_{k=1}^{n-1}
\RR_{n-k}(y+\xi,y) (f_k(z-\xi)-f_k(-z-\xi))\, d\xi.
\end{split}
\ee
If $\beta_u$ exists, where $3<u\le 4$, then 
$$
\int_0^\infty (t+1) f(\xi+t)\, dt = \PP \{ X_1 \ge \xi \} +
\int_0^\infty \PP \{ X_1 \ge \xi+t \}\, dt \ll \frac{1}{(\xi+1)^{u-1}}.
$$
Therefore, by Theorem \ref{xlarge},
\be\label{intf}
\RR_n(y+\xi,y) \ll \frac{y+1}{n^{3/2}(1+\xi)^{u-1}}.
\ee
Let $V_1$ be the contribution to the integral in \eqref{Rmain} from 
$1\le k\le n_0$, let $V_2$ be the contribution from $n_0+1 \le k \le
n/2$ and $V_3$ is the contribution from $n/2 < k \le n-1$.
By \eqref{intf},
\be\label{S1}
V_1 \ll \frac{y+1}{n^{3/2}}
\sum_{k=1}^{n_0} \int_0^\infty f_k(z-\xi)+f_k(-z-\xi) \,
d\xi = 2n_0 \frac{y+1}{n^{3/2}}.
\ee
When $k\ge n_0+1$, Lemma \ref{fn1} implies that
\be\label{fk2}
\begin{split}
f_k(z-&\xi)-f_k(-z-\xi) = \frac{e^{-\frac{1}{2k}(z-\xi)^2}}{\sqrt{2\pi
    k}} \( 1 - e^{-2\xi z/k} \) + O\pfrac{1}{k^{(u-1)/2}}\\
&\qquad + O\biggl[ \(\frac{|z-\xi|}{k^{3/2}} + \frac{|z-\xi|^3}{k^{5/2}} \) 
  e^{-(z-\xi)^2/2k} 
+ \( \frac{z+\xi}{k^{3/2}} + \frac{(z+\xi)^3}{k^{5/2}} \)
    e^{-(z+\xi)^2/2k} \biggr] \\
&\ll \frac{1}{k^{(u-1)/2}} + \frac{(z+1)(\xi+1)}{k^{3/2}} e^{-(z-\xi)^2/2k}. 
\end{split}
\ee
By \eqref{intf}, we have
\begin{align*}
V_2 &\ll \frac{y+1}{n^{3/2}} \sum_{n_0+1 \le k \le n/2}
\int_0^\infty \frac{1}{k^{(u-1)/2}(\xi+1)^{u-1}} + \frac{z+1}{k^{3/2}(\xi+1)^{u-2}}
e^{-(z-\xi)^2/2k} \, d\xi \\
&\ll \frac{y+1}{n^{3/2}} \left[ 1 +
  (z+1) \sum_{k\le n/2} \frac{1}{k^{3/2}} \int_0^\infty
  \frac{1}{(\xi+1)^{u-2}} e^{-(z-\xi)^2/2k} \, d\xi \right].
\end{align*}
The integral on the right side is
\begin{align*}
&\le e^{-z^2/8k} \int_0^{z/2} \frac{d\xi}{(\xi+1)^{u-2}} +
\int_{-z/2}^\infty \frac{e^{-w^2/2k}}{(z+w)^{u-1}}\, dw \\
&\ll  e^{-z^2/8k} + \min\( \frac{1}{(z+1)^{u-3}},
\frac{k^{1/2}}{(z+1)^{u-2}} \).
\end{align*}
Hence
\be\label{S2}
\begin{split}
V_2 &\ll \frac{y+1}{n^{3/2}} (z+1) \sum_{k=1}^\infty
  k^{-3/2} \( e^{-z^2/8k} +  \min \( \frac{1}{(z+1)^{u-3}},
  \frac{k^{1/2}}{(z+1)^{u-2}} \) \) \\
&\ll \frac{y+1}{n^{3/2}} (z+1) \left[ \frac{1}{z+1} +
    \frac{1}{(z+1)^{u-2}} \sum_{k\le z^2} \frac{1}{k} \right] 
\ll \frac{y+1}{n^{3/2}}. 
\end{split}
\ee
By Lemma \ref{jumpsize} and \eqref{fk2},
\be\label{S3}
V_3 \ll \( \frac{z+1}{n^{3/2}} + \frac{1}{n^{(u-1)/2}} \)
\sum_{j=1}^\infty \int_0^\infty (\xi+1)
\RR_j(y+\xi,y)\, d\xi \ll \frac{z+1}{n^{3/2}} + \frac{1}{n^{(u-1)/2}}.
\ee
Putting together \eqref{Rmain}, \eqref{intf}, \eqref{S1}, \eqref{S2}
and \eqref{S3}, we arrive at
$$
\RR_n(x,y) = f_n(x) - f_n(y+z) + O\( \frac{y+z+1}{n^{3/2}} +
\frac{1}{n^{(u-1)/2}}\).
$$
Since $|x| \le M \sqrt{n}$, Lemma \ref{fn1} implies
$f_n(x) \gg n^{-1/2}$ for sufficiently large $n$, 
the implied constant depending on the
distribution of $X_1$ and also on $M$.  Hence
$$
R_n(x,y) = 1 - \frac{f_n(y+z)}{f_n(x)} + O\( \frac{y+z+1}{n} +
\frac{1}{n^{(u-2)/2}}\).
$$
Finally, by Lemma \ref{fn1} again,
\begin{align*}
\frac{f_n(y+z)}{f_n(x)} &= e^{-\frac{1}{2n}((y+z)^2-x^2)} +
  O\( \frac{|x| + y + z + 1}{n} + \frac{1}{n^{(u-2)/2}} \) \\
&= e^{-2yz/n} + O\( \frac{y+z+1}{n} + \frac{1}{n^{(u-2)/2}} \).
\end{align*}
Again the implied constant depends on $M$.
This concludes the proof of Theorem \ref{posz}.

\bigskip
{\bf Acknowledements.}
The author thanks Valery Nevzorov for suggesting to utilize the 
reflection principle in a form similar to that in Lemma \ref{recur2}.
The author is grateful to the referees for carefully reading the paper
and for several small corrections and suggestions.


\bibliographystyle{amsplain}
\bibliography{condprob}

\end{document}